\documentclass[11pt]{article}
\usepackage{epsfig}

\newtheorem{Theorem}{Theorem}[section]

\newcommand{\me}{{-1}}

\newcommand{\iy}{\infty}

\newcommand{\al}{\alpha}

\newcommand{\eps}{\varepsilon}
\newcommand{\ga}{\gamma}

\newcommand{\de}{\delta}
\newcommand{\la}{\lambda}
\newcommand{\ph}{\varphi}
\newcommand{\vro}{\varrho}
\newcommand{\tht}{\theta}
\newcommand{\si}{\sigma}
\newcommand{\om}{\omega}

\newcommand{\vsk}{\vspace{2mm}}
\newcommand{\vsg}{\vspace{3mm}}

\begin{document}

{\Large\bf
\begin{center}
On the Eigenvalues of Certain Canonical Higher-Order
Ordinary Differential Operators
\end{center}}

{ \bf
\begin{center}
Albrecht B\"ottcher and Harold Widom
\end{center}}

\begin{quote}
\renewcommand{\baselinestretch}{1.0}
\footnotesize
We consider the operator of taking the $2p$th derivative of a function
with zero boundary conditions for the function and its first $p-1$ derivatives
at two distinct points.
Our main result provides an asymptotic formula for the eigenvalues and
resolves a question on the appearance of certain regular numbers in the
eigenvalue sequences for $p=1$ and $p=3$.
\end{quote}


\vsk

\section{{\large Introduction}}\label{S1}

\noindent
Let $\al$ be a natural number. We consider the eigenvalue problem
\begin{eqnarray}
& & \hspace{-10mm} (-1)^\al u^{(2\al)}(x)= \la \,u(x) \quad \mbox{for}
\quad x \in [0,1], \label{1.1}\\
& & \hspace{-10mm} u(0)=u'(0)=\ldots=u^{(\al-1)}(0)=0, \quad
u(1)=u'(1)=\ldots=u^{(\al-1)}(1)=0. \label{1.2}
\end{eqnarray}
This problem has countably many eigenvalues, which are all positive and
converge to infinity.
We denote the sequence of the eigenvalues by
$\{\lambda_{n,\alpha}\}_{n=n_0}^\infty$ where $n_0$ will be
chosen in dependence on $\alpha$ (the a priori choice
$n_0=1$ will turn out to be inconvenient). Thus, the first
eigenvalue is $\lambda_{n_0,\alpha}$, the second is
$\lambda_{n_0+1,\alpha}$, and so on. We also
put $\mu_{n,\alpha}
=\sqrt[2\al]{\la_{n,\al}}$.
In \cite{BoWi2} we determined the asymptotics of the minimal eigenvalue
$\la_{\min, \al}$ as $\al$ goes to infinity. The result is
\begin{equation}
\la_{\min,\al} = \sqrt{8\pi\al}\,\left(\frac{4\al}{e}\right)^{2\al}
\left(1+O\left(\frac{1}{\sqrt{\al}}\right)\right).\label{asy}
\end{equation}
In \cite{BoWi2} we also observed that $\la_{\min,3}$ is very close to
$(2\pi)^6$. We here show that $\la_{\min,3}$ is in fact equal to
$(2\pi)^6$. It is furthermore well known that $\la_{\min,1}=\pi^2$.
The purpose of this paper is to give an answer to the
question whether these coincidences are accidents or not.
We shall prove that they are due to the accidents that $\cos\pi=-1$
and $\cos\frac{\pi}{3}=\frac{1}{2}$ are nonzero rational numbers.

\vsk
For $\al=1$, the eigenvalues of (\ref{1.1}), (\ref{1.2}) are known
to be $\la_{n,1}=(n\pi)^2$ ($n \ge 1$). If $\al=2$, the eigenvalues
are given by $\la_{n,2}=\mu_{n,2}^4$ where $\{\mu_{n,2}\}$ is the sequence
of the positive solutions of the equation
\begin{equation}
\cos \mu=\frac{1}{\cosh \mu}. \label{1.3a}
\end{equation}
This was shown in \cite{Hor} and \cite{Par99}. From (\ref{1.3a}) we infer that
if $n_0$ is appropriately chosen, then the sequence $\{\mu_{n,2}\}_{n=n_0}^\iy$
has the asymptotics $\mu_{n,2}=\frac{\pi}{2} +n\pi+\de_n$ with $\de_n
\sim 2\,(-1)^{n+1}e^{-\pi/2}e^{-n\pi}$ as $n \to \iy$. As usual, $x_n \sim y_n$
means that $x_n/y_n \to 1$ as $n \to \iy$.
Thus, in contrast to the case $\al=1$, $\{\mu_{n,2}\}_{n=n_0}^\iy$ does not contain
an arithmetic progression. In both \cite{Hor} and \cite{Par99} it was also established
that $\mu_{\min,2}=4.7300$ and, accordingly, $\la_{\min,2}=500.5467$. We may therefore take
$n_0=1$ and write the eigenvalue sequence in the form
$\{(\pi/2+n\pi)^4+\xi_n)\}_{n=1}^\iy$ with $\xi_n \sim 4\pi^3n^3\de_n$.

\vsk
Now let $\al=3$. The
general solution of equation (\ref{1.1}) is
\[u(x)=\sum_{j=0}^5 C_j\exp(\mu \eps^{2j+1}x) \quad \mbox{with}
\quad \eps=e^{\pi i/6}, \quad \mu=\sqrt[6]{\la}.\]
Consequently, the boundary conditions (\ref{1.2}) are satisfied if and only if
the determinant of the matrix
\begin{equation}
A_3(\mu)=\left(\begin{array}{llllll}
1 & 1 & 1 & 1 & 1 & 1 \\
\eps & \eps^3 & \eps^5 & \eps^7 & \eps^9 & \eps^{11}\\
\eps^2 & \eps^6 & \eps^{10} & \eps^{14} & \eps^{18} & \eps^{22}\\
a\om & \om^2 & b\om & b\om^{-1} & \om^{-2} & a \om^{-1}\\
a\om\eps & \om^2\eps^3 & b\om\eps^5 & b\om^{-1}\eps^7 & \om^{-2}\eps^9 & a \om^{-1}\eps^{11}\\
a\om\eps^2 & \om^2\eps^6 & b\om\eps^{10} & b\om^{-1}\eps^{14}
& \om^{-2}\eps^{18} & a \om^{-1}\eps^{22}\end{array}\right)\label{1.4}
\end{equation}
is zero, where
\begin{equation}
a=e^{\mu\,{\rm Re}\,\eps}=e^{\mu\,\sqrt{3}/2},\quad
b=e^{-\mu\,{\rm Re}\,\eps}=e^{-\mu\,\sqrt{3}/2},\quad
\om=e^{\mu\,{\rm Im}\,\eps}=e^{\mu\pi i/2}. \label{1.3}
\end{equation}

Figure \ref{miplo} shows the minimum of the absolute values of
the eigenvalues of the matrix $A_3(\mu)$
in dependence on $\mu$.
The unit on the horizontal axis is $\pi$. We see that the minimum of the moduli of the
eigenvalues and hence the determinant of $A_3(\mu)$ is zero for $\mu_{n,3}$
sharply concentrated at
the values of $n\pi$ with $n \ge 2$. The question is whether $\mu_{n,3}$
is exactly $n\pi$ or not. Theorem \ref{Th 1.1} provides the answer.

{\vspace*{7cm}
\begin{figure}[ht]
\begin{picture}(4,4)
\centerline{{\epsfig{file=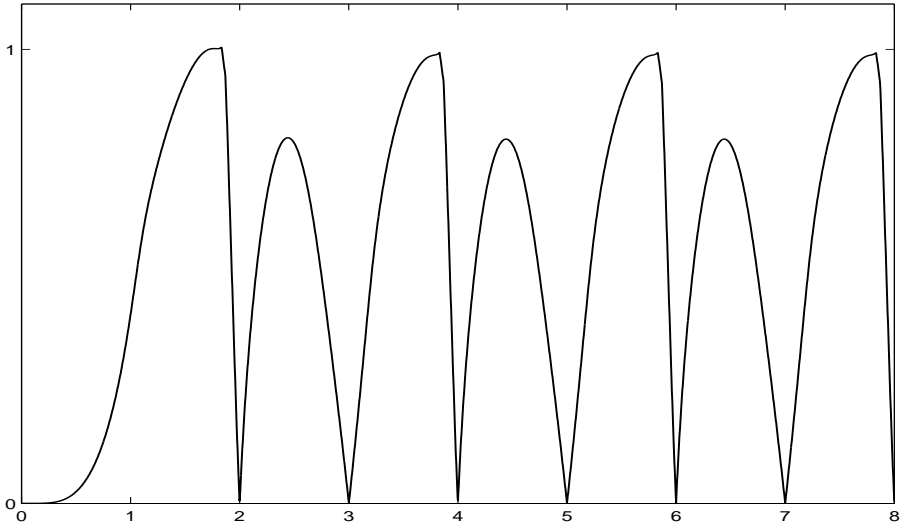, height=7cm, width=12cm}}}
\end{picture}
\caption{The minimum of the moduli of the eigenvalues of $A_3(\mu)$ in dependence of $\mu$.}
\label{miplo}
\end{figure}}

\begin{Theorem} \label{Th 1.1}
For $\al=3$, the eigenvalues of (\ref{1.1}), (\ref{1.2}) are
$\{\la_{n,3}\}_{n=2}^\iy$ with
\[\la_{n,3}=(n\pi)^6 \quad \mbox{if}\quad  n \ge 2 \quad \mbox{is even}\]
and
\[\la_{n,3}=(n \pi +\de_n)^6 \quad \mbox{if} \quad n \ge 3 \quad \mbox{is odd},\]
where the $\de_n$'s are nonzero numbers satisfying
$\de_n \sim 8\,(-1)^{\lfloor n/2 \rfloor+1}e^{-(\pi\sqrt{3}/2)n}$.
Here $\lfloor n/2 \rfloor$ stands for the
integral part of $n/2$. In particular,
$\mu_{n,3}=n\pi$ if and only if $n$ is even.
\end{Theorem}

Mark Embree computed the first five roots of equation
(\ref{2.1}) and thus the numbers $\mu_{n,3}$ for $n=3,5,7,9,11$ up to ten correct digits
after the comma. The result is in the table. The last column
shows the values of $8\,(-1)^{\lfloor n/2 \rfloor+1}e^{-(\pi\sqrt{3}/2)n}$.

\begin{center}
\begin{tabular}{|l|rr|c|}
\hline
\rule[-1.5ex]{0ex}{4.5ex}
& $\mu_{n,3}$ & & $8\,(-1)^{\lfloor n/2 \rfloor+1}e^{-(\pi\sqrt{3}/2)n}$\\
\hline
$n=3$ & $9.4270555708$ & $=3\,\pi+0.0022776101$ & $+\,0.0022821082$\\
$n=5$ & $15.7079533785$ & $=5\,\pi-0.0000098894$ & $-\,0.0000098893$\\
$n=7$ & $21.9911486179$ & $=7\,\pi+0.0000000428$ & $+\,0.0000000428$\\
$n=9$ & $28.2743338821$ & $=9\,\pi-0.0000000002$ & $-\,0.0000000002$\\
$n=11$ & $34.5575191894$ & $=11\,\pi+0.0000000000$ & $+\,0.0000000000$\\
\hline
\end{tabular}
\end{center}

If $\al=1$ then $\{\mu_{n,1}\}_{n=1}^\iy$ is an arithmetic progression, and if
$\al=3$ then $\{\mu_{n,3}\}_{n=2}^\iy$ contains an arithmetic progression.
The following result involves the reason for these two peculiarities.
Actually it does more. It gives the asymptotics of the eigenvalues
of problem (\ref{1.1}), (\ref{1.2}) for arbitrary $\al$.

\begin{Theorem} \label{Th 1.2}
If $\al \ge 5$ is odd, one can choose
$n_0$ so that the sequence $\{\mu_{n,\al}\}_{n=n_0}^\iy$
satisfies
\begin{equation}
\mu_{n,\al}=n\pi+\frac{2\,(-1)^n}{\sin^2\!\frac{\pi}{2\al}}\,
e^{-n\pi\sin\!\frac{\pi}{\al}}\,\sin\left(n\pi\cos\frac{\pi}{\al}\right)
+O\left(e^{-n\pi\sin\!\frac{2\pi}{\al}}\right).
\label{1.5}
\end{equation}
If $\al \ge 4$ is even, there is an $n_0$ such that $\{\mu_{n,\al}\}_{n=n_0}^\iy$
can be written as
\begin{equation}
\mu_{n,\al}=\frac{\pi}{2}+n\pi+\frac{2\,(-1)^{n+1}}{\sin^2\!\frac{\pi}{2\al}}\,
e^{-\left(\frac{\pi}{2}+n\pi\right)\sin\!\frac{\pi}{\al}}\,
\cos\left(\left(\frac{\pi}{2}+n\pi\right)\cos\frac{\pi}{\al}\right)
+O\left(e^{-2\,n\pi\sin\!\frac{\pi}{\al}}\right).
\label{1.6}
\end{equation}
The sequence $\{\mu_{n,\al}\}_{n=n_0}^\iy$ contains an arithmetic progression
if and only if $\cos\frac{\pi}{\al}$ is a nonzero rational number, that is,
if and only if $\al=1$ or $\al=3$.
\end{Theorem}

We remark that (\ref{1.5}) is also true for $\al=3$, but that in this case it
amounts to $\mu_{n,3}=n\pi+O\left(e^{-(\pi \sqrt{3}/2)n}\right)$, which is weaker than
Theorem \ref{Th 1.1}. Formula (\ref{1.6}) becomes valid for $\al=2$ after replacing
the numerator $2\,(-1)^{n+1}$ by $(-1)^{n+1}$. The reason for this discrepancy will be given
at the end of Section \ref{S3}.

\vsk
In the cases $\al=1,2,3$ the values of $\mu_{\min,\al}$ are $\pi$, $4.7300$
($\approx \frac{3}{2}\,\pi$), $2\pi$, respectively. This suggests that the
$n_0$ in Theorem \ref{Th 1.2} is the integral part of $\frac{\al+1}{2}$ and
hence asymptotically equals $\frac{\al}{2}$ as $\al \to \iy$. However, (\ref{asy})
implies that $\mu_{\min,\al} \sim \frac{4\al}{e}$ and thus $n_0\pi \sim
\frac{4\al}{e}$, which shows that actually $n_0$ is asymptotically equal
to $\frac{4}{\pi e}\,\al
=0.4684\,\al$. It was just these tiny but significant differences that have
kindled our interest in the subject in \cite{BoWi2} and here.

\vsk
The properties of problem (\ref{1.1}), (\ref{1.2})
depend on whether $\al$ is
odd or even. We therefore study these two cases separately.

\section{{\large The odd case}}\label{S2}

There remains nothing to say on the case $\al=1$. So let us begin with $\al=3$.
The following theorem provides us with a formula for the determinant of matrix
(\ref{1.4}) for arbitrary $a,b,\om$,
that is, for $a,b,\om$ that are not necessarily of the form (\ref{1.3}).

\begin{Theorem} \label{Th 2.1}
Let $\eps=e^{\pi i/6}$ and let
$a$,$b$, and $\om\neq 0$ be arbitrary complex numbers. Then the determinant of matrix
(\ref{1.4}) is
\begin{eqnarray*}
{\rm det}\,A_3(\mu) & = & 12\,ab\,(a+b)(\om-\om^{-1}) +3\,(a^2+b^2)(\om^{-2}-\om^2)\\
& & +3\,ab\,(\om^4-8\,\om^2+8\,\om^{-2}-\om^{-4})+12\,(a+b)(\om-\om^{-1})\\
& = & 12\,ab\,(a+b)(\om-\om^{-1}) -3\,(a^2+b^2)(\om-\om^{-1})(\om+\om^{-1})\\
& & +3\,ab\,(\om -\om^{-1})(\om+\om^{-1})(\om^2+\om^{-2}-8)+12\,(a+b)(\om-\om^{-1}).
\end{eqnarray*}
\end{Theorem}

\noindent
{\em Proof.}
Let $V_{ijk}$ denote the determinant of the (Vandermonde) matrix that is
constituted by the first three
rows and the columns $i,j,k$ of $A_3(\mu)$.
We expand the determinant of $A_3(\mu)$ by its last three rows
using Laplace's theorem and group the $\left({6 \atop 3}\right)=20$ terms so that
we get a polynomial in $a$ and $b$. What results is
\begin{eqnarray*}
& & a^2b\,(-\om \,V_{136}V_{245}+\om^{-1}\,V_{146}V_{235})\\
& & +ab^2\,(-\om\,V_{134}V_{256}+\om^{-1}\,V_{346}V_{125})\\
& & +a^2\,(\om^2\,V_{126}V_{345}-\om^{-2}\,V_{156}V_{234})\\
& & +b^2\,(\om^2\,V_{234}V_{156}-\om^{-2}\,V_{345}V_{126})\\
& & +ab\, (-\om^4\,V_{123}V_{456} +\om^2\,V_{124}V_{356} -\om^{-2}\,V_{145}V_{236}
+ \om^2\,V_{236}V_{145}\\
& & \quad -V_{246}V_{135}+V_{135}V_{246}-\om^{-2}\,V_{356}V_{124}
+\om^{-4}\,V_{456}V_{123})\\
& & +a\,(-\om\,V_{125}V_{346}+\om^{-1}\,V_{256}V_{134})\\
& & +b\,(-\om\,V_{235}V_{146}+\om^{-1}\,V_{245}V_{136}).
\end{eqnarray*}
The determinants $V_{ijk}$ are
\begin{eqnarray*}
& & V_{136}=V_{235}=V_{145}=2\sqrt{3}\,\eps^2,\quad
V_{245}=V_{146}=V_{236}=2\sqrt{3}\,\eps^4,\\
& & V_{256}=V_{346}=V_{124}=2\sqrt{3}\,\eps^8,\quad
V_{134}=V_{125}=V_{356}=2\sqrt{3}\,\eps^{10},\\
& & V_{126}=V_{234}=V_{456}=\sqrt{3},\quad
V_{345}=V_{156}=V_{123}=-\sqrt{3},\\
& & V_{246}=-\frac{3}{2}\,\sqrt{3}, \quad V_{135}=\frac{3}{2}\,\sqrt{3}.
\end{eqnarray*}
Inserting these values in the above expression for the determinant
of $A_3(\mu)$, we arrive
at the asserted formula. $\; \rule{2mm}{2mm}$

\begin{Theorem} \label{Th 2.2}
For $\al=3$, the sequence of the eigenvalues of (\ref{1.1}), (\ref{1.2})
is $\{\la_{n,3}\}_{n=2}^\iy =\{\mu_{n,3}^6\}_{n=2}^\iy$ where $\mu_{n,3}=n\pi$
if $n$ is even and $\mu_{3,3}, \mu_{5,3}, \mu_{7,3}, \ldots$ are the positive
solutions of the equation
\begin{equation}
\cos\frac{\mu}{2}=\frac{4\,\cosh\frac{\mu\sqrt{3}}{2}}{\cosh\mu\sqrt{3}}
+\frac{1}{\cosh\mu\sqrt{3}}\,\left[\cos\frac{\mu}{2}\,\cos\mu -4 \,\cos\frac{\mu}{2}
\right]. \label{2.1}
\end{equation}
\end{Theorem}

\noindent
{\em Proof.} We apply Theorem \ref{Th 2.1} to the case where $a,b,\om$ are
given by (\ref{1.3}). Since in that case $ab=1$ and
\begin{eqnarray*}
& & a+b=2\,\cosh\frac{\mu\sqrt{3}}{2}, \quad a^2+b^2=2\,\cosh\mu\sqrt{3},\\
& & \om-\om^{-1}=2i\,\sin\frac{\mu}{2}, \quad \om+\om^{-1}=2\,\cos\frac{\mu}{2},
\quad \om^2+\om^{-2}=2\,\cos \mu,
\end{eqnarray*}
we obtain that
\begin{equation}
{\rm det}\,A_3(\mu)=24\,i\,\sin\frac{\mu}{2}\,\left[
4\,\cosh\frac{\mu\sqrt{3}}{2}-\cosh\mu\sqrt{3}\,\cos\frac{\mu}{2}
+\cos\frac{\mu}{2}\,\cos\mu-4\cos\frac{\mu}{2}
\right].\label{2.2}
\end{equation}
The factor $\sin\frac{\mu}{2}$ produces the eigenvalues
$\la_{2k,3}=(2k\pi)^6$ ($k=1,2,3, \ldots$).
The term in the brackets of (\ref{2.2}) is
$\cosh\mu\sqrt{3}$ times
\begin{equation}
\frac{4\,\cosh\frac{\mu\sqrt{3}}{2}}{\cosh\mu\sqrt{3}}
+\frac{1}{\cosh\mu\sqrt{3}}\,\left[\cos\frac{\mu}{2}\,\cos\mu -4 \,\cos\frac{\mu}{2}
\right]-\cos\frac{\mu}{2}, \label{2.3}
\end{equation}
and the zeros of (\ref{2.3}) are the solutions of (\ref{2.1}). $\; \rule{2mm}{2mm}$

\begin{figure}[ht]
\begin{picture}(4,4)
\put(30,0){{\epsfig{file=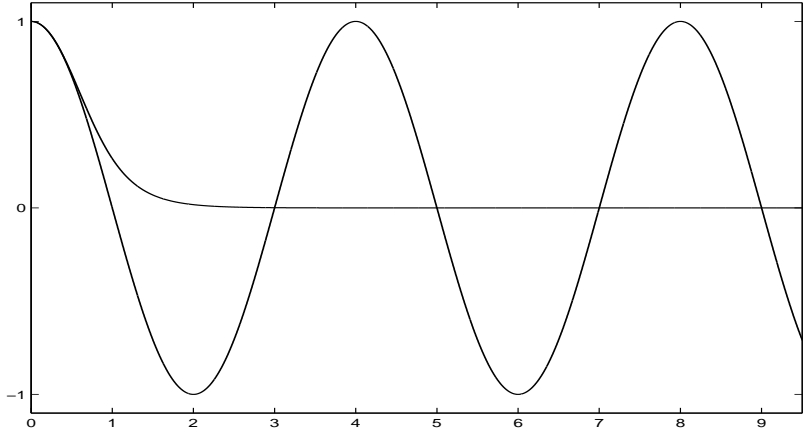, height=6cm, width=12cm}}}
\end{picture}
\caption{The graphs of the functions on the left and right of (\ref{2.1}),
the unit on the horizontal axis being $\pi$.}
\label{gra}
\end{figure}

\vsk
Theorem \ref{Th 1.1} is almost straightforward from Theorem \ref{Th 2.2}.
Figure \ref{gra} shows the graphs of the functions on the left and right
of (\ref{2.1}). We see that the first intersection occurs at approximately
$3\pi$.
Thus, the zeros of (\ref{2.3})
may be written as
\[\mu_{n,3}=n\pi+\de_n \quad (n=3,5,7, \ldots)\]
with $\de_n\to 0$
as $n \to \iy$. The right-hand side of (\ref{2.1}) has the asymptotics
$4\,e^{-(\sqrt{3}\pi/2)\mu}$. It follows that
\[4e^{-(\pi\sqrt{3}/2)n} \sim \cos\frac{\mu_{n,3}}{2}
=\cos\frac{n\pi+\de_n}{2}=-\sin\frac{n\pi}{2}\sin\frac{\de_n}{2}\sim
(-1)^{\lfloor n/2\rfloor +1}\frac{\de_n}{2}\]
and, in particular, that $\de_n \neq 0$ for all sufficiently
large $n$. As (\ref{2.3}) is definitely nonzero for $\mu=3\pi, 5\pi,
7\pi, \ldots$, we conclude that $\de_n \neq 0$ for {\em all} $n$.

\pagebreak

\vspace*{5cm}

\vsk
We now turn to general odd numbers $\al \ge 5$. One can build up the analogue
$A_\al(\mu)$ of matrix (\ref{1.4}). The following theorem gives the first two
terms of the asymptotics of ${\rm det}\,A_\al(\mu)$ as $\mu \to \iy$.

\begin{Theorem} \label{Th 2.3}
If $\al =2k+1 \ge 5$, then
\begin{equation}
{\rm det}\,A_\al(\mu)=K_1e^{\ga_1\mu}\sin\mu+K_2e^{\ga_2\mu}
\sin(s_{2k-1}\mu)+O\left(e^{\ga_3\mu}\right) \label{2.4}
\end{equation}
where $K_1, K_2$ are nonzero constants satisfying $-K_2/K_1=2/\sin^2\!\frac{\pi}{2\al}$,
\begin{eqnarray*}
& & \ga_1=\si+2c_{2k-3}+2c_{2k-1}, \quad
\ga_2=\si+2c_{2k-3}+c_{2k-1},\quad \ga_3=\si+c_{2k-3}+2c_{2k-1},\\
& & \si=2\sum_{j=1}^{k-2} c_{2j-1}, \quad c_\ell=\cos\frac{\ell\pi}{2\al},
\quad s_\ell=\sin\frac{\ell\pi}{2\al}
\quad (\ell=1,3, \ldots,2k-1).
\end{eqnarray*}
\end{Theorem}

\noindent
{\em Proof.} To avoid heavy notation and to show the essence of the matter,
we restrict ourselves to the case where $\al=5$. The general solution of (\ref{1.1})
is
\[u(x)=\sum_{j=0}^9C_j\exp\left(\mu\,\eps^{2j+1}x\right),
\quad \eps=e^{\pi i/10}, \quad \mu=\sqrt[10]{\la}.\]
Thus, the analogue $A_5(\mu)$ of matrix (\ref{1.4}) is a $10 \times 10$ matrix
of a structure completely analogous to the one of (\ref{1.4}). The first $5$
rows are of Vandermonde type, the second row being
\[\left(\; \eps \quad \eps^3 \quad \eps^5 \quad \eps^7 \quad \eps^9
\quad \eps^{11} \quad \eps^{13}
\quad \eps^{15} \quad \eps^{17} \quad \eps^{19}\;\right).\]
The remaining $5$ rows result from the first $5$ rows by multiplying each column
by a factor. These factors appear without the powers of $\eps$ in the $6$th row,
which is
\[\left(\;
a\om \quad b\tau \quad \tht \quad b^\me\tau \quad a^\me\om
\quad a^\me\om^\me \quad b^\me\tau^\me
\quad \tht^\me \quad b\tau^\me \quad a\om^\me\;\right)
\]
with
\begin{eqnarray*}
& & a=e^{\mu c_1}, \quad b=e^{\mu c_3}, \quad \om=e^{i\mu s_1},
\quad \tau=e^{i\mu s_3}, \quad \tht=e^{i \mu},\\
& & c_1=\cos\frac{\pi}{10}, \quad c_3=\cos\frac{3\pi}{10}, \quad
s_1=\sin\frac{\pi}{10}, \quad s_3=\sin\frac{3\pi}{10}.
\end{eqnarray*}
To make things absolutely safe, we still note that the first and
second columns
of $A_5(\mu)$ are
\begin{eqnarray*}
& & \left(\begin{array}{llllllllll}
1 & \eps & \eps^2 & \eps^3 & \eps^4 &
a\om  & a\om \eps & a\om \eps^2 & a\om \eps^3 & a\om\eps^4\end{array}\right)^\top,\\
& & \left(\begin{array}{llllllllll}
1 & \eps^3 & \eps^6 & \eps^{9} & \eps^{12} &
b\tau  & b\tau \eps^3 & b\tau \eps^6 & b\tau \eps^9 & b\tau\eps^{12}\end{array}\right)^\top.
\end{eqnarray*}
We expand ${\rm det}\,A_5(\mu)$ by its last $5$ rows using the Laplace theorem.
Taking into account that
\[a^2b^2=e^{\mu(2c_1+2c_3)}, \quad a^2b=e^{\mu(2c_1+c_3)}, \quad
a^2=e^{\mu\cdot 2c_1}, \quad ab^2=e^{\mu(c_1+2c_3)}\]
and that $c_1+2c_3>2c_1$, we see that the asymptotics of ${\rm det}\,A_5(\mu)$
is
\[L_1\,a^2b^2 +L_2\, a^2b+O(ab^2).\]
We are left with determining $L_1$ and $L_2$. Let $V_{j_1, \ldots, j_5}$ denote
the determinant of the Vandermonde matrix at the intersection of the first
$5$ rows and the columns $j_1, \ldots, j_5$ of $A_5(\mu)$. With $X:=10$, we have
\begin{eqnarray*}
L_1 & = & -V_{1239X}V_{45678}\,\tht+V_{1289X}V_{34567}\,\tht^\me,\\
L_2 & = & V_{1238X}V_{45679}\,\tau + V_{1249X}V_{35678}\,\tau
- V_{1279X}V_{34568}\,\tau^\me -V_{1389X}V_{24567}\,\tau^\me.
\end{eqnarray*}
It remains to compute the products of the form $V_MV_{M'}$ where $M'=\{1,2, \ldots, 10\}
\setminus M$. Obviously,
\begin{eqnarray*}
V_M & = & \prod_{{j,\ell \in M \atop j<\ell}}\left(\eps^{2\ell-1}-\eps^{2j-1}\right)
=\prod_{{j,\ell \in M \atop j<\ell}}\left(\eps^{2\ell-1}+\eps^{2j-1+10}\right)\\
& = & \prod_{{j,\ell \in M \atop j<\ell}}\vro_{\ell-j}\,\eps^{(2\ell-1+2j-1+10)/2}
= \prod_{{j,\ell \in M \atop j<\ell}}\vro_{\ell-j}\,\eps^{\ell+j+4}
\end{eqnarray*}
where
\[\vro_{\ell-j}=\left|\eps^{2\ell-1}-\eps^{2j-1}\right|=2\sin\frac{(\ell-j)\pi}{10}.\]
Letting
\[R_M=\prod_{{j,\ell \in M \atop j<\ell}}\vro_{\ell-j}\,
\prod_{{j,\ell \in M' \atop j<\ell}} \vro_{\ell-j}\]
we get $V_MV_{M'}=R_M\eps^q$ with
\[q=4\,(1+2+\ldots+10)+2\left({5 \atop 2}\right)\cdot 4 =300.\]
Since $\eps^{20}=1$, it follows that $\eps^{300}=1$ and thus
$V_MV_{M'}=R_M$. A direct computation shows that
\begin{eqnarray*}
& & R_{1239X}= R_{1289X} = -\vro_1^8\,\vro_2^6\,\vro_3^4\,\vro_4^2 =:-S,\\
& & R_{1238X} = R_{1249X} =R_{1279X} = R_{1389X}
= \vro_1^6\,\vro_2^6\,\vro_3^4\,\vro_4^2\,\vro_5^2 =:T,
\end{eqnarray*}
whence
\[
L_1  = S(\tht-\tht^\me)=2iS\,\sin\mu,\quad
L_2  =  2T(\tau-\tau^\me)=4iT\,\sin(\mu s_3).
\]
In summary,
\[
{\rm det}\,A_5(\mu)=2iS\,e^{\mu(2c_1+2c_3)}\sin\mu
+4iT\,e^{\mu(2c_1+c_3)}\sin(\mu s_3)+O\left(e^{\mu(c_1+2c_3)}\right),
\]
and since $-4iT/2iS=2\vro_5^2/\vro_1^2=2/\sin^2\!\frac{\pi}{10}$, this is
(\ref{2.4}) for $\al=5$. $\; \rule{2mm}{2mm}$

\vsk
We are now in a position to prove Theorem \ref{Th 1.2} for odd $\al$.
By Theorem \ref{Th 2.3}, the equation
${\rm det}\,A_\al(\mu)=0$ is of the form
\begin{equation}
\sin \mu =-\frac{K_2}{K_1}\,e^{-(\ga_1-\ga_2)\mu}\,\sin(s_{2k-1}\mu)
+O\left(e^{-(\ga_1-\ga_3)\mu}\right). \label{2.6}
\end{equation}
The solutions of (\ref{2.6}) are $n\pi+\de_n$ with small $\de_n$'s satisfying
\[\sin(n\pi+\de_n)=(-1)^n \sin\de_n=O\left(e^{-n\pi(\ga_1-\ga_2)}\right),\]
which implies that $\de_n=
O\left(e^{-n\pi(\ga_1-\ga_2)}\right)$. Consequently,
\begin{eqnarray*}
& & (-1)^n\left(\de_n+O\left(e^{-2\,n\pi(\ga_1-\ga_2)}\right)\right)\\
& & \quad
=-\frac{K_2}{K_1}\,e^{-(n\pi+\de_n)(\ga_1-\ga_2)}\,\sin(s_{2k-1}n\pi+s_{2k-1}\de_n)
+O\left(e^{-n\pi(\ga_1-\ga_3)}\right),
\end{eqnarray*}
and since $e^{-\de_n(\ga_1-\ga_2)}=1+O\left(e^{-n\pi(\ga_1-\ga_2)}\right)$ and
\begin{eqnarray*}
& & \sin (s_{2k-1}n\pi+s_{2k-1}\de_n)\\
& & \quad
=\sin(s_{2k-1}n\pi)\,\cos(s_{2k-1}\de_n)+\cos(s_{2k-1}n\pi)\,\sin(s_{2k-1}\de_n)\\
& & \quad
=\sin(s_{2k-1}n\pi)+O\left(e^{-n\pi(\ga_1-\ga_2)}\right),
\end{eqnarray*}
we arrive at the formula
\begin{equation}
(-1)^n\de_n=-\frac{K_2}{K_1}\,e^{-n\pi(\ga_1-\ga_2)}\,\sin(s_{2k-1}n\pi)
+O\left(e^{-n\pi \ga}\right) \label{n14}
\end{equation}
with $\ga=\min(\ga_1\!-\!\ga_3, 2(\ga_1\!-\!\ga_2))$. Because
\[2(\ga_1\!-\!\ga_2)=2c_{2k-1}=2\,\sin\frac{\pi}{\al} >
\sin\frac{2\pi}{\al} =c_{2k-3}=\ga_1\!-\!\ga_3,\]
we see that
$\ga=\ga_1\!-\!\ga_3=\sin\frac{2\pi}{\al}$. Finally, taking into account that
$-K_2/K_1=2/\sin^2\!\frac{\pi}{2\al}$ and $s_{2k-1}=\cos\frac{\pi}{\al}$,
we obtain (\ref{1.5}) from (\ref{n14}).

\vsk
Now suppose $\{\mu_{n,\al}\}$ contains an arithmetic progression.
Because $\mu_{n,\al}=n\pi+o(1)$ by virtue
of (\ref{1.5}), this progression must be of the form $\{n\ell\pi\}_{n=N}^\iy$
with some natural number $\ell$. As $\sin(n\ell\pi)=0$, we obtain from the
second order asymptotics of Theorem \ref{Th 2.3} that
\[K_2e^{\ga_2n\ell\pi}\sin(s_{2k-1}n\ell\pi)+o\left(e^{\ga_2n\ell \pi}\right)=0\]
for all $n \ge N$, which implies that
$\sin(s_{2k-1}n\ell\pi)\to 0$ as $n \to \iy$.
Consequently, $s_{2k-1}n\ell$ must converge to zero modulo $1$. By Kronecker's theorem,
this is only possible
if $s_{2k-1}\ell$ and thus $s_{2k-1}$ is a rational number. Since
$s_{2k-1}=\cos\frac{\pi}{\al}$,
we conclude that $\cos\frac{\pi}{\al}$ has to be rational. In an appendix
we prove that this is not the case for $\al \ge 4$.

\section{{\large The even case}}\label{S3}

Let $\al \ge 4$ be an even number. Here is the analogue of Theorem \ref{Th 2.3}.

\begin{Theorem} \label{Th 3.1}
If $\al=2k \ge 4$, then
\[{\rm det}\,A_\al(\mu)=K_1e^{\ga_1\mu}\cos\mu +K_2e^{\ga_2\mu}\cos(s_{k-1}\mu)
+O\left(e^{\ga_3 \mu}\right)\]
where $K_1, K_2$ are nonzero constants such that
$-K_2/K_1=2/\sin^2\!\frac{\pi}{2\al}$,
\begin{eqnarray*}
& & \ga_1=\si+2c_{k-1}, \quad
\ga_2=\si+c_{k-1}, \quad \ga_3=\si,\\
& & \si=1+2\sum_{j=1}^{k-2}c_j, \quad c_j=\cos\frac{j\pi}{\al},
\quad s_j=\sin\frac{j\pi}{\al}
\quad (j=1,2, \ldots,k-1).
\end{eqnarray*}
\end{Theorem}

\noindent
{\em Proof.} We confine ourselves to the case $\al=4$. The matrix $A_4(\mu)$ is
\[\left(\begin{array}{llllllll}
1 & 1 & 1 & 1 & 1 & 1 & 1 & 1\\
1 & \eps & \eps^2 & \eps^3 & \eps^4 & \eps^5 & \eps^6 & \eps^7\\
1 & \eps^2 & \eps^4 & \eps^6 & \eps^8 & \eps^{10} & \eps^{12}& \eps^{14}\\
1 & \eps^4 & \eps^6 & \eps^9 & \eps^{12} & \eps^{15} & \eps^{18}& \eps^{21}\\
a & b\om & \tau & b^\me\om & a^\me & b^\me\om^\me & \tau^\me & b\om^\me\\
a & b\om\eps & \tau\eps^2 & b^\me\om\eps^3 & a^\me\eps^4
& b^\me\om^\me\eps^5 & \tau^\me\eps^6 & b\om^\me\eps^7\\
a & b\om\eps^2 & \tau\eps^4 & b^\me\om\eps^6 & a^\me\eps^8
& b^\me\om^\me\eps^{10} & \tau^\me\eps^{12} & b\om^\me\eps^{14}\\
a & b\om\eps^4 & \tau\eps^6 & b^\me\om\eps^9 & a^\me\eps^{12}
& b^\me\om^\me\eps^{15} & \tau^\me\eps^{18} & b\om^\me\eps^{21}
\end{array}\right)\]
with
\begin{eqnarray*}
& & \eps=e^{\pi i/4}, \quad a=e^\mu, \quad b=e^{\mu c_1}, \quad
\om=e^{i\mu s_1}, \quad \tau=e^{i\mu},\\
& & c_1=\cos\frac{\pi}{4}, \quad s_1=\sin\frac{\pi}{4}.
\end{eqnarray*}
Laplace expansion through the last $4$ rows yields
\[{\rm det}\,A_4(\mu)=L_1\,ab^2+L_2\,ab +O(a)\]
with
\begin{eqnarray}
L_1 & = & V_{1238}V_{4567}\,\tau + V_{1278}V_{3456}\,\tau^\me,\nonumber\\
L_2 & = & -V_{1237}V_{4568}\,\om -V_{1248}V_{3567}\,\om
-V_{1268}V_{3457}\,\om^\me -V_{1378}V_{2456}\,\om^\me. \label{four}
\end{eqnarray}
We have
\[V_MV_{M'}=\eps^6\,\prod_{{j,\ell \in M \atop j<\ell}}\vro_{\ell-j}
\prod_{{j,\ell \in M' \atop j<\ell}}\vro_{\ell-j} =:\eps^6\,R_M=-iR_M\]
where
\[\vro_{\ell-j}=\left|\eps^{\ell-1}-\eps^{j-1}\right|=2\sin\frac{(\ell-j)\pi}{8}.\]
This gives
\[
R_{1238} =R_{1278} =\vro_1^6\,\vro_2^4\,\vro_3^2 =:S,\quad
R_{1237} =R_{1248} =R_{1268} =R_{1378}=
\vro_1^4\,\vro_2^4\,\vro_3^2\,\vro_4^2=:T.
\]
We finally obtain that
\begin{eqnarray*}
{\rm det}\,A_4(\mu) & = & -iS\,ab^2(\tau+\tau^\me)+2iT\,ab(\om+\om^\me) +O(a)\\
& = & -2iS\,e^{\mu(1+2c_1)}\cos\mu+4iT\,e^{\mu(1+c_1)}\cos(s_1\mu)+O(e^\mu)
\end{eqnarray*}
with $-4iT/(-2iS)=2\vro_4^2/\vro_1^2=2/\sin^2\!\frac{\pi}{8}$. $\; \rule{2mm}{2mm}$

\vsk
Armed with Theorem \ref{Th 3.1} we can prove Theorem \ref{Th 1.2}
for even numbers $\al$. The equation ${{\rm det}}\,A_\al(\mu)=0$ reads
\[\cos\mu=-\frac{K_2}{K_1}\,e^{-(\ga_1-\ga_2)\mu}\,\cos(s_{k-1}\mu)+
O\left(e^{-(\ga_1-\ga_3)\mu}\right)\]
and we have $\ga_1\!-\!\ga_2=c_{k-1}$ and $\ga_1\!-\!\ga_3=2c_{k-1}$.
The solutions of this equation are $\frac{\pi}{2}+n\pi+\de_n$ with
$\de_n=O\left(e^{-n\pi c_{k-1}}\right)$. It follows that
\begin{eqnarray*}
& & (-1)^{n+1}\sin\de_n=\cos\left(\frac{\pi}{2}+n\pi+\de_n\right)\\
& & \quad
=-\frac{K_2}{K_1}\,e^{-\left(\left(\frac{\pi}{2}+n\pi+\de_n\right)c_{k-1}\right)}\,
\cos\left(s_{k-1}\left(\frac{\pi}{2}+n\pi+\de_n\right)\right)+
O\left(e^{-2n\pi c_{k-1}}\right),
\end{eqnarray*}
and since
\begin{eqnarray*}
& & \sin\de_n=\de_n+O\left(e^{-2n\pi c_{k-1}}\right),\quad
e^{-\de_nc_{k-1}}=1+O\left(e^{-n\pi c_{k-1}}\right),\\
& & \cos\left(s_{k-1}\left(\frac{\pi}{2}+n\pi+\de_n\right)\right)
=\cos\left(s_{k-1}\left(\frac{\pi}{2}+n\pi\right)\right)
+O\left(e^{-n\pi c_{k-1}}\right),
\end{eqnarray*}
we get the representation
\[(-1)^{n+1}\de_n=-\frac{K_2}{K_1}
\,e^{-\left(\left(\frac{\pi}{2}+n\pi\right)c_{k-1}\right)}\,
\cos\left(s_{k-1}\left(\frac{\pi}{2}+n\pi\right)\right)+
O\left(e^{-2n\pi c_{k-1}}\right).\]
It remains to notice that $-K_2/K_1=2/\sin^2\!\frac{\pi}{2\al}$,
$s_{k-1}=\cos\frac{\pi}{\al}$, and $c_{k-1}=\sin\frac{\pi}{\al}$.

\vsk
If $\{\mu_{n,\al}\}$ contains an arithmetic
progression, we can argue as in Section \ref{S2} to see that there is a natural
number $\ell$ such that
\[\frac{1}{2}\,s_{k-1}+\ell ns_{k-1} \to \frac{1}{2}\quad \mbox{modulo} \;1\]
as $n \to \iy$. Kronecker's theorem implies again that $s_{k-1}=\cos\frac{\pi}{\al}$
must be rational.

\vsk
Finally, in the case $\al=2$ we have to deal with the matrix
\[A_2(\mu)=\left(\begin{array}{llll}
1 & 1 & 1 & 1\\1 & \eps & \eps^2 & \eps^3\\
a & \tau & a^\me & \tau^\me\\
a & \tau\eps & a^\me \eps^2 & \tau^\me \eps^3\end{array}\right),\]
where $\eps = e^{\pi i/2}$, $a=e^\mu$, $\tau=e^{\mu i}$. The determinant
of this matrix is $L_1 a+L_2 +O(a^\me)$ with
\[L_1=V_{12}V_{34}\tau+ V_{34}V_{12}\tau^\me, \quad L_2=-V_{13}V_{24}
-V_{24}V_{13}.\]
We see that the constant $L_2$ is the sum of two terms, which is in contrast
to the case $\al \ge 4$, where the constant $L_2$ is the sum of four terms
as in (\ref{four}). This explains why for $\al=2$ the numbers $\mu_{n,2}$
are $\frac{\pi}{2}+n\pi$ plus half of the subsequent term of
(\ref{1.6}). Incidentally, a straightforward computation yields
${\rm det}\,A_2(\mu)=8(1-\cos \mu \cosh \mu)$, which leads to
(\ref{1.3a}).

\vspace{5mm}
\noindent
{\bf \large Appendix}

\vsg
\noindent
Here is, just for completeness, a proof of the fact that $\cos\frac{\pi}{\al}$
is irrational for $\al \ge 4$. Since $\cos(n\tht)$ is a polynomial with integer coefficients
of $\cos\tht$, the rationality of $\cos\tht$ implies that of $\cos(n\tht)$.
We are therefore left with proving that $\cos\frac{\pi}{\al}$ is irrational
if $\al$ is $4,6,9$ or a prime number $p \ge 5$.

\vsk
The cases $\al=4$ and $\al=6$ are trivial.
So let $\al$ be $9$ or a prime number $p \ge 5$ and put $x=\cos\frac{\pi}{\al}$.
We denote by $\ph$ the polynomial $\ph(y)=2y^2-1$ and by $\ph^n$ the $n$th iterate
of $\ph$, $\ph^2(y)=\ph(\ph(y))$ and so on. Clearly $\ph^n(x)=\cos\frac{2^n\pi}{\al}$.
If $n=6$ then $2^n=1$ modulo $9$, and if $n=p-1$ then $2^n=1$ modulo $p$ (Fermat).
Thus, for these $n$ we have $\ph^n(x)=\cos\frac{(\ell \al+1)\pi}{\al}=\pm x$. Since
$\ph^n(0)=1$ for $n \ge 2$ and the leading coefficient $a_m$ of $\ph^n(y)$ is a power
of $2$, $a_m=2^k$, it follows that $x$ satisfies an algebraic equation of the
form $2^kx^m + \ldots +1=0$. The only rational solutions of such an equation
are of the form $x=\pm 1/2^j$. Because $1>\cos\frac{\pi}{\al} \ge
\cos\frac{\pi}{5} >\frac{1}{2}$, we arrive at the conclusion that $\cos\frac{\pi}{\al}$
must be irrational.

\vsg
\noindent
\begin{minipage}[t]{8cm}
Albrecht B\"ottcher\\
Fakult\"at f\"ur Mathematik\\
TU Chemnitz\\
09107 Chemnitz\\
Germany\\[1ex]
aboettch@mathematik.tu-chemnitz.de
\end{minipage}
\begin{minipage}[t]{6.5cm}
Harold Widom \\
Department of Mathematics\\
University of California\\
Santa Cruz, CA 95064\\
USA\\[1ex]
widom@math.ucsc.edu

\end{minipage}

\end{document}